\documentclass[12pt]{amsart}
\usepackage{amsfonts,amssymb}
\usepackage[mathscr]{eucal}
\usepackage{amsmath, amsthm}
\usepackage{mathrsfs}
\usepackage{amsbsy}
\usepackage{dsfont}
\usepackage{bbm}
\makeindex
\input xypic
\xyoption{all}
\begin{document}
\baselineskip = 5mm
\newcommand \lra {\longrightarrow}
\newcommand \hra {\hookrightarrow}
\newcommand \ZZ {{\mathbb Z}} 
\newcommand \NN {{\mathbb N}} 
\newcommand \QQ {{\mathbb Q}} 
\newcommand \RR {{\mathbb R}} 
\newcommand \CC {{\mathbb C}} 
\newcommand \sgr {{\mathfrak S}} 
\newcommand \bcA {{\mathscr A}}
\newcommand \bcB {{\mathscr B}}
\newcommand \bcC {{\mathscr C}}
\newcommand \bcD {{\mathscr D}}
\newcommand \bcE {{\mathscr E}}
\newcommand \bcF {{\mathscr F}}
\newcommand \bcU {{\mathscr U}}
\newcommand \C {{\mathscr C}}
\newcommand \X {{\mathscr X}}
\newcommand \im {{\rm im}}
\newcommand \Hom {{\rm Hom}}
\newcommand \colim {{{\rm colim}\, }} 
\newcommand \iHom {{\underline {\rm Hom}}} 
\newcommand \End {{\rm {End}}}
\newcommand \Ob {{\rm {Ob}}}
\newcommand \coker {{\rm {coker}}}
\newcommand \id {{\rm {id}}}
\newcommand \cm {{\mathscr M}}
\newcommand \uno {{\mathbbm 1}} 
\newcommand \Le {{\mathbbm L}} 
\newcommand \PR {{\mathbb P}} 
\newcommand \Spec {{\rm {Spec}}}
\newcommand \SP {{SP}} 
\newcommand \Pic {{\rm {Pic}}}
\newcommand \Alb {{\rm {Alb}}}
\newcommand \Corr {{Corr}}
\newcommand \Sym {{\rm {Sym}}}
\newcommand \cha {{\rm {char}}}
\newcommand \tr {{\rm {tr}}} 
\newcommand \res {{\rm {res}}} 
\newtheorem{theorem}{Theorem}
\newtheorem{lemma}[theorem]{Lemma}
\newtheorem{sublemma}[theorem]{Sublemma}
\newtheorem{corollary}[theorem]{Corollary}
\newtheorem{example}[theorem]{Example}
\newtheorem{exercise}[theorem]{Exersize}
\newtheorem{proposition}[theorem]{Proposition}
\newtheorem{remark}[theorem]{Remark}
\newtheorem{notation}[theorem]{Notation}
\newtheorem{definition}[theorem]{Definition}
\newtheorem{conjecture}[theorem]{Conjecture}
\newenvironment{pf}{\par\noindent{\em Proof}.}{\hfill\framebox(6,6)
\par\medskip}
\title[Codimension two algebraic cycles on threefolds over a field]
{\bf On the continuous part of codimension two algebraic cycles on
threefolds over a field}
\author{V. Guletski\u \i }

\date{16 July 2007}

\begin{abstract}
\noindent Let $X$ be a non-singular projective threefold over an
algebraically closed field of any characteristic, and let $A^2(X)$
be the group of algebraically trivial codimension $2$ algebraic
cycles on $X$ modulo rational equivalence with coefficients in $\QQ
$. Assume $X$ is birationally equivalent to another threefold $X'$
admitting a fibration over an integral curve $C$ whose generic fiber
$X'_{\bar \eta }$, where $\bar \eta =\Spec (\overline {k(C)})$,
satisfies the following three conditions: (i) the motive $M(X'_{\bar
\eta })$ is finite-dimensional, (ii) $H^1_{et}(X_{\bar \eta },\QQ
_l)=0$ and (iii) $H^2_{et}(X_{\bar \eta },\QQ _l(1))$ is spanned by
divisors on $X_{\bar \eta }$. We prove that, provided these three
assumptions, the group $A^2(X)$ is representable in the weak sense:
there exists a curve $Y$ and a correspondence $z$ on $Y\times X$,
such that $z$ induces an epimorphism $A^1(Y)\to A^2(X)$, where
$A^1(Y)$ is isomorphic to $\Pic ^0(Y)$ tensored with $\QQ $. In
particular, the result holds for threefolds birational to
three-dimensional Del Pezzo fibrations over a curve.
\end{abstract}





\maketitle

\section{Introduction}
\label{intro}

Let $X$ be a non-singular projective variety of dimension $d$ over
an algebraically closed field $k$. For any integer $0\leq i\leq d$
let $CH^i(X)$ be a Chow group of codimension $i$ algebraic cycles on
$X$ with coefficients in $\QQ $ modulo rational equivalence
relation. Let $A^i(X)$ be a subgroup in $CH^i(X)$ generated by
cycles algebraically equivalent to zero. Recall that an algebraic
cycle $Z$ on $X$ is called to be algebraically equivalent to zero if
it can be deformed to $0$ in a trivial family over a non-singular
projective curve. In codimension one $A^1(X)$ can be identified with
the Picard group $\Pic ^0(X)$ tensored with $\QQ $, and $A^d(X)$ is
the group of zero-dimensional algebraic cycles of degree zero modulo
rational equivalence relation. The group $A^i(X)$ is said to be {\it
(weakly) representable} if there exists a curve $Y$ and a cycle
class $z$ in $CH^i(Y\times X)$, such that the induced homomorphism
$z_*:A^1(Y)\to A^i(X)$ is surjective. For example, $A^1(X)$ is
representable if $X$ is a curve. For a surface, if $k$ is
uncountable, representability of $A^2(X)$ is equivalent to the
triviality of the kernel $T(X)$ of the Albanese homomorphism
$A^2(X)\to \Alb (X)\otimes \QQ $. If $k=\CC $, representability of
zero cycles on a surface without non-trivial globally holomorphic
$2$-forms was conjectured by S.Bloch, and it was a part of the
intuition leading to the whole Bloch-Beilinson motivic vision of
algebraic cycles, \cite{Jannsen2}. On the other hand, if $X$ is a
surface with $p_g>0$, the group $A^2(X)$ is too far from to be
representable, \cite{Mumford}. Recently it was discovered that
$A^2(X)$ is representable for a surface $X$ with $p_g=0$ if and only
if the motive $M(X)$ of $X$ (with coefficients in $\QQ $) is
finite-dimensional in the sense of S.Kimura, \cite{GP2},
\cite{Kimura}. The aim of the present paper is to show that motivic
finite-dimensionality can be also useful for the study of $A^2(X)$
in three-dimensional case. We will show that for a certain type of
threefolds $X$ representability of $A^2(X)$ follows from
finite-dimensionality of the motive of the generic fiber of an
appropriate fibration of $X$ over a curve. In full generality this
statement should not be true, of course.

To state our result precisely we need to fix some notation. Let
$f:X\to C$ be a fibering of a threefold over an integral curve,
$\eta =\Spec (k(C))$ be the generic point on $C$, $\bar \eta =\Spec
(\overline {k(C)})$ the spectrum of the closure of $k(C)$, and
$X_{\bar \eta }$ the generic fiber of the morphism $f$ over $\bar
\eta $. Let also $H^1_{et}(X_{\bar \eta },\QQ _l)$ and
$H^2_{et}(X_{\bar \eta },\QQ _l(1))$ be the $i$-th etale $l$-adic
cohomology groups, where $l$ is a prime different from the
characteristic of the ground field $k$.

\begin{theorem}
\label{main} Let $X$ be a non-singular projective threefold over an
algebraically closed field $k$ of any characteristic. Assume that
$X$ is birationally equivalent, by means of a chain of blowing ups
and blowing downs with non-singular centers, to another threefold
$X'$ fibered over an integral curve both defined over $k$. Assume,
furthermore, that the generic fiber $X'_{\bar \eta }$ of that
fibering satisfies the following three assumptions:

\begin{itemize}

\item[{\it (i)}]
the motive $M(X'_{\bar \eta })$ is finite-dimensional;

\item[{\it (ii)}]
$H^1_{et}(X'_{\bar \eta },\QQ _l)=0$ and

\item[{\it (iii)}]
$H^2_{et}(X'_{\bar \eta },\QQ _l(1))$ is spanned by divisors on
$X'_{\bar \eta }$.

\end{itemize}

\noindent Then the group $A^2(X)$ is representable.
\end{theorem}

Let us indicate how much that theorem is applicable. First of all,
we claim the result over an arbitrary ground field $k$, i.e. we do
not make any restriction on its characteristic. If $k=\CC $ then,
according to the Minimal Model Program over $\CC $, all non-singular
projective threefolds can be divided into two parts: those which are
birational to a $\QQ $-factorial $X'$ with at most terminal
singularities and $K_{X'}$ nef, and those which are birational to a
$\QQ $-factorial $X'$ with at most terminal singularities with an
extremal contraction $X'\to C$ being a Del Pezzo fibration. The
second case can be divided into $3$ subcases: if $\dim (C)=0$ then
$X'$ is Fano, if $\dim (C)=2$ then $X'$ is a conic bundle, and if
$\dim (C)=1$ then the generic fiber $X'_{\eta }$ is a non-singular
projective Del Pezzo surface over $\eta $. In the last case $X$
satisfies the assumptions of Theorem \ref{main}, whence
representability of $A^2(X)$. A typical example here is a general
non-singular cubic $X$ in $\PR ^4$. The group $A^2(X)$ for $X$ is
known to be representable (H.Clemens, P.Griffiths) and the motive of
the generic hyperplane section is finite-dimensional. If we take a
general quartic $X\subset \PR ^4$, which is Fano, then $A^2(X)$ is
representable as well (S.Bloch, J.Murre), but its non-singular
hyperplane sections are K3's - those surfaces whose motivic
finite-dimensionality is not known yet.

The paper is divided into two parts. In Section \ref{prelim} we
recall some known facts and also prove preliminary results on
representability of the continuous part of codimension two algebraic
cycles on threefolds. Section \ref{proofmain} is the main part of
the paper where we make main computations and prove Theorem
\ref{main}.

\section{\it Preliminary results}
\label{prelim}

Let $S$ be a non-singular connected quasi-projective variety over a
field $k$, and let $\SP (S)$ be the category of all non-singular and
projective schemes over $S$. Assume we are given two objects of that
category, and let $X=\cup _jX_j$ be the connected components of $X$.
For any non-negative $m$ let
  $$
  \Corr _S^m(X,Y)=\oplus _jCH^{e_j+m}(X_j\times _SY)
  $$
be the group of relative correspondences of degree $m$ from $X$ to
$Y$ over $S$, where $e_j$ is the relative dimension of $X_j$ over
$S$. For example, given a morphism $f:X\to Y$ in $\SP (S)$, the
transpose $\Gamma _f^t$ of its graph $\Gamma _f$ is in $\Corr
_S^0(X,Y)$. For any two correspondences $f:X\to Y$ and $g:Y\to Z$
their composition $g\circ f$ is defined, as usual, by the formula
  $$
  g\circ f={p_{13}}_*(p_{12}^*(f)\cdot p_{23}^*(g))\; ,
  $$
where the central dot denotes the intersection of cycle classes in
the sense of \cite{Ful}. The category $\cm (S)$ of Chow-motives over
$S$ with coefficients in $\QQ $ can be defined then as a
pseudoabelian envelope of the category of correspondences with
certain ``Tate twists" indexed by integers. For any non-singular
projective $X$ over $S$ its motive $M(X/S)$ is defined by the
relative diagonal $\Delta _{X/S}$, and for any morphism $f:X\to Y$
in $\SP (S)$ the correspondence $\Gamma _f^{t}$ defines a morphism
$M(f):M(Y/S)\to M(X/S)$. The category $\cm (S)$ is rigid with a
tensor product satisfying the formula
  $$
  M(X/S)\otimes M(Y/S)=M(X\times _SY)\; ,
  $$
so that the functor
  $$
  M:\SP (S)^{\rm op}\lra \cm (S)
  $$
is tensor. The scheme $S/S$ indexed by $0$ gives the unite $\uno _S$
in $\cm (S)$, and when it is indexed by $-1$, it gives the Lefschetz
motive $\Le _S$. If $E$ is a multisection of degree $w>0$ of $X/S$,
we set
  $$
  \pi _0=\frac {1}{w}\, [E\times _SX]\; \; \; \hbox{and}\; \; \;
  \pi _{2e}=\frac {1}{w}\, [X\times _SE]\; ,
  $$
where $e$ is the relative dimension of $X/S$. Then one has the
standard isomorphisms $\uno _S\cong (X,\pi _0,0)$ and $\Le
_S^{\otimes e}\cong (X,\pi _{2e},0)$. Finally, if $f:T\to S$ is a
morphism of base schemes over $k$, then $f$ gives a base change
tensor functor $f^*:\cm (S)\to \cm (T)$. All the details about Chow
motives over a non-singular base can be found, for instance, in
\cite{GM}.

Below we will use basic facts from the theory of finite dimensional
motives, or, more generally, finite dimensional objects, see
\cite{Kimura} or \cite{AK}. Roughly speaking, once we have a tensor
$\QQ $-linear pseudo-abelian category $\bcC $, one can define wedge
and symmetric powers of any object in $\bcC $. Then we say that
$X\in \Ob (\bcC )$ is finite-dimensional, \cite{Kimura}, if it can
be decomposed into a direct sum, $X=Y\oplus Z$, such that $\wedge
^mY=0$ and $\Sym ^nZ=0$ for some non-negative integers $m$ and $n$.
The property to be finite dimensional is closed under direct sums
and tensor products, etc. A morphism $f:X\to Y$ in $\bcC $ is said
to be numerically trivial if for any morphism $g:Y\to X$ the trace
of the composition $g\circ f$ is equal to zero,
\cite[7.1.1-7.1.2]{AK}. The important role in the below arguments is
played by the following result:

\begin{proposition}
\label{AKnilp} Let $\bcC $ be a tensor $\QQ $-linear pseudo-abelian
category, let $X$ be a finite-dimensional object in $\bcC $, and let
$f$ be a numerically trivial endomorphism of $X$. Then $f$ is
nilpotent in the ring $\End (X)$.
\end{proposition}

\begin{pf}
See \cite[7.5]{Kimura} for Chow motives and \cite[9.1.14]{AK} in the
abstract setting.
\end{pf}

\begin{lemma}
\label{traces} Let $F:\bcC _1\to \bcC _2$ be a tensor functor
between two rigid tensor $\QQ $-linear pseudo-abelian categories.
Assume $F$ induces an injection $\End (\uno _{\bcC
_1})\hookrightarrow \End (\uno _{\bcC _2})$. Then, if $X$ is a
finite-dimensional object in $\bcC _1$ and $F(X)=0$, it follows that
$X=0$ as well.
\end{lemma}

\begin{pf}
Let $g$ be an endomorphism of $X$. Then $F(\tr (g\circ 1_X))=\tr
(F(g\circ 1_X))$, see \cite{DeligneMilne}, page 116. Since $F(X)=0$,
we have $F(g\circ 1_X)=0$. Then $\tr (g\circ 1_X)=0$ because $F$ is
an injection on the rings of endomorphisms of units. Hence, the
identity morphism $1_X$ is numerically trivial. Since $X$ is finite
dimensional, it is trivial by Proposition \ref{AKnilp}.
\end{pf}

\medskip

The following lemma was pointed out to me by J. Ayoub:

\begin{lemma}
\label{limit} Let $S$ be a non-singular integral variety over $k$,
let $\eta $ be the generic point of $S$, and let $W$ be the set of
Zariski open subsets in $S$. Then
  $$
   \cm (\eta )=\colim _{U\in W}\cm (U)\; ,
  $$
where the colimit is taken in the $2$-category of pseudo-abelian
tensor categories with coefficients in $\QQ $.
\end{lemma}

\begin{pf}
Assume we are given with a tensor $\QQ $-linear pseudo-abelian
category $\mathscr T$ and a set of tensor functors
  $$
  G_U:\cm (U)\lra \mathscr T\; ,
  $$
compatible with restriction functors
  $$
  i^*:\cm (U)\lra \cm (V)
  $$
for each inclusion of Zariski open subsets $i:V\subset U$. We need
to show that there exists a unique functor
  $$
  F:\cm (\eta )\lra {\mathscr T}\; ,
  $$
such that the composition of a pull-back $\cm (U)\to \cm (\eta )$
with $F$ coincides with $G_U$ for each $U$ from $W$. The proof is,
actually, just a systematic use of spreads of algebraic cycles and
the localization sequence for Chow groups. Indeed, let $M=(X,p,n)$
be a motive over $\eta $, and let $X'$ and $p'$ be spreads of $X$
and $p$ respectively over some Zariski open subset $U$ in $S$.
Shrinking $U$ if necessary and applying the localization for Chow
groups we may assume that $p'^2=p'$, i.e. $p'$ is a relative
projector over $U$. The triple $M'=(X',p',n)$ is an object in $\cm
(U)$, so that we define $F(M)$ to be $G_U(M')$. Systematically
shrinking Zariski open subsets one can easily show that such defined
$F(M)$ does not depend on the choice of spreads. On morphisms $F$
can be defined in a similar way, because morphisms in Chow motives
are just algebraic cycles. The uniqueness is evident when taking
into account the localization for Chow groups again.
\end{pf}

\medskip

Also we need to recall an equivalent reformulation of
representability of the continuous part in $A^2$ for threefolds. Let
$X$ be a non-singular projective threefold $X$ over an algebraically
closed field $k$. As we pointed out in Introduction, the group
$A^2(X)$ is representable if there exists a non-singular projective
curve $Y$ and a correspondence $a\in CH^2(Y\times X)$, such the
homomorphism $a_*:A^1(Y)\to A^2(X)$ induced by $a$ is surjective,
where $A^1(Y)$ is $\Pic ^0(Y)\cong \Alb (Y)$ tensored with $\QQ $.
By technical reasons it is convenient to introduce another one
definition: $A^2(X)$ is representable if there exist a finite
collection of non-singular projective curves $Y_1,\dots ,Y_m$ and
correspondences $a_i\in CH^2(Y_i\times X)$, such that the
homomorphism
  $$
  \sum _{i=1}^m(a_i)_*:\oplus _{i=1}^mA^1(Y_i)\lra A^2(X)
  $$
is surjective. Evidently, the first definition implies the second
one. The following argument for the inverse implication has been
taken from \cite{Bloch}. By Bertini's theorem, we take a
non-singular one-dimensional linear section $Y$ of $Y_1\times \dots
\times Y_m$. The inclusion $Y\hra Y_1\times \dots \times Y_m$ gives
rise to a surjection on Albanese varieties $\Alb (Y)\to \Alb
(Y_1\times \dots \times Y_m)$. Fix a closed point $y_i$ on each
$Y_i$. For any $1\leq i\leq m$ the points $y_1,\dots
,y_{i-1},y_{i+1},\dots ,y_m$ give the embedding $Y_i\hra Y_1\times
\dots \times Y_m$ sending any point $y\in Y_i$ into the point
$(y_1,\dots ,y_{i-1},y,y_{i+1},\dots ,y_m)$, where $y$ stays on
$i$-th place. Composing such an embedding with the corresponding
projection we get the identity map on $Y_i$. Applying the functor
$\Alb $ we get the homomorphisms $\Alb (Y_i)\to \Alb (Y_1\times
\dots \times Y_m)$ induced by the above embeddings and the
homomorphisms $\Alb (Y_1\times \dots \times Y_m)\to \Alb (Y_i)$
induced by the projections. Since a finite product of abelian groups
is also a coproduct, we get two homomorphisms $\Alb (Y_1)\times
\dots \times \Alb (Y_m)\to \Alb (Y_1\times \dots \times Y_m)$ and
$\Alb (Y_1\times \dots \times Y_m)\to \Alb (Y_1)\times \dots \times
\Alb (Y_m)$ whose composition is identity on $\Alb (Y_1)\times \dots
\times \Alb (Y_1)$. It follows that the second homomorphism is
surjective. Composing it with the above surjective homomorphism
$\Alb (Y)\to \Alb (Y_1\times \dots \times Y_m)$ we get a surjective
homomorphism $\Alb (Y)\to \times _{i=1}^m\Alb (Y_i)=\oplus
_{i=1}^m\Alb (Y_i)$. Since $Y$ and all $Y_1,\dots ,Y_m$ are
non-singular projective curves, we can replace $\Alb $ by $A^1$
getting a surjective homomorphism $A^1(Y)\to \oplus
_{i=1}^mA^1(Y_i)$. Composing it with the above homomorphism $\sum
_{i=1}^m(a_i)_*:\oplus _{i=1}^mA^1(Y_i)\to A^2(X)$ we get a
surjective homomorphism $A^1(Y)\to A^2(X)$, as required. Thus, both
definitions are equivalent.

\section{\it The proof of Theorem \ref{main}}
\label{proofmain}

Let $X$ be a non-singular projective threeefold over an
algebraically closed field $k$ satisfying the assumptions of Theorem
\ref{main}. The below proof can be divided into two parts - chain of
reductions and computations with relative correspondences. The
reductions consist, roughly speaking, in finite extensions of the
base curve $C$, cutting finite collections of points out of $C$ and,
respectively, removing their fibers out of $X$. Then, of course, we
will work with a new family dealing with $A^2$ in the non-compact
case. The point is that representability of $A^2(X)$ for $X\to C$ is
local on the base.

To be more precise, we start with the following two lemmas both
proved in \cite{Bloch}:

\begin{lemma}
\label{finite} Let $Y\to X$ be a morphism of non-singular projective
threefolds of finite degree over $k$. If $A^2(Y)$ is representable,
then $A^2(X)$ is representable as well.
\end{lemma}

\begin{lemma}
\label{birat} Let $X$ and $X'$ be two non-singular projective
threefolds over a field. Assume that $X$ is birationally equivalent
to $X'$ by means of a chain of blowing ups and blowing downs with
non-singular centers. Then representability of $A^2(X)$ is
equivalent to representability of $X$.
\end{lemma}

\bigskip

\noindent {\it Reduction 1}

\medskip

By Lemma \ref{birat} we can replace $X$ by $X'$ from the very
beginning (see the formulation of Theorem \ref{main}). In
particular, $X$ is equipped now with the fibering $X\to C$ over a
non-singular projective curve $C$, whose generic fiber satisfies the
three assumptions in Theorem \ref{main}.

\bigskip

\noindent {\it Reduction 2}

\medskip

Let $\eta =\Spec (k(C))$ be the generic point of the map
  $$
  f:X\lra C\; ,
  $$
and let $\bar \eta $ be the spectrum of an algebraic closure
$\overline {k(C)}$ of the function field $k(C)$. We assume that
$H^1_{et}(X_{\bar \eta },\QQ _l)=0$ and that $H^2_{et}(X_{\bar \eta
},\QQ _l(1))$ is algebraic. Let $b_2$ be the second Betti number for
the generic fiber $X_{\bar \eta }$ and let
  $$
  D_1,\dots ,D_{b_2}
  $$
be divisors over $\bar \eta $ generating $H^2_{et}(X_{\bar \eta
},\QQ _l)$. We can extend $C$ by a finite map $C'\to C$ so that, if
$X'=X\times _CC'$ is a base change, the divisors $D_i$ are defined
over the generic point $\eta '$ of the curve $C'$. Since $X'\to X$
is a morphism of finite degree, we may assume that the divisors
$D_1,\dots ,D_{b_2}$ are defined over $\eta $ by Lemma \ref{finite}.

\bigskip

Now we need the following easy lemma:

\begin{lemma}
\label{uptofin} Let $X$ be a non-singular projective threefold over
$k$ and let $A^2(X)=V\oplus W$ be a splitting of the $\QQ $-vector
space $A^2(X)$ into two subspaces. Assume, furthermore, that $V$ is
finite-dimensional. Then $A^2(X)$ is representable if and only if
there exists a finite collection of non-singular projective curves
$Y_1,\dots ,Y_m$ and correspondences $a_i\in CH^2(Y_i\times X)$,
such that the homomorphism $\sum _{i=1}^m(a_i)_*:\oplus
_{i=1}^mA^1(Y_i)\to A^2(X)$ is onto the subspace $W$. In other
words, representability of $A^2(X)$ is up to a finite-dimensional
subspace in $A^2(X)$.
\end{lemma}

\begin{pf}
Let $\{ v_1,\dots ,v_n\} $ be a basis for $V$. For each index $j$
let $B_j$ be an algebraically equivalent to zero algebraic cycle
representing the cycle class $v_j$. By the definition of algebraic
equivalence there exists a non-singular projective curve $T_j$, two
closed points $p_j$ and $q_j$ on $T_j$, and an algebraic cycle $Z_j$
on $T_j\times X$, such that $Z_j$ intersects the divisors $p_j\times
X$ and $q_j\times X$ properly, and these intersections are $B_j$ and
$0$ respectively. In that case, if $z_j$ is a class of the cycle
$A_j$, the value of the homomorphism $(z_j)_*$ on the class of the
algebraically trivial zero-cycle $p_j-q_j$ is exactly the cycle
class $v_j$. Thus, the homomorphism $\sum (z_j)_*:\oplus A^1(T_j)\to
A^2(X)$ covers the space $V$. To complete the proof we need only to
add the curves $Y_i$ from the assumptions and use the second
equivalent definition of representability of $A^2(X)$.
\end{pf}

\bigskip

\noindent {\it Reduction 3}

\medskip

Let $\{ p_1,\dots ,p_m\} $ be a finite set of closed points on the
curve $C$, let
  $$
  U=C-\{ p_1,\dots ,p_m\}
  $$
be the non-singular locus of the family, and let
  $$
  Y=X-f^{-1}(U)\; .
  $$
By the localization exact sequence
  $$
  \oplus _{j=1}^mCH^1(X_{p_j})\lra CH^2(X)\lra CH^2(Y)\to 0
  $$
the $\QQ $-vector space $CH^2(X)$ splits as
  $$
  CH^2(X)=CH^2(Y)\oplus I\; ,
  $$
where
  $$
  I=\im (\oplus _jCH^1(X_{p_j})\lra CH^2(X))\; .
  $$
In other words, the localization homomorphism $CH^2(X)\to CH^2(Y)$
has a section $CH^2(Y)\to CH^2(X)$. It is not hard to see that for
any cycle class $\alpha \in A^2(Y)$ one can find a cycle class
$\beta $ in the preimage of $\alpha $ with respect to the surjective
homomorphism $CH^2(X)\lra CH^2(Y)$, such that $\beta $ is
algebraically trivial as well. Then we have a surjective
localization homomorphism
  $$
  A^2(X)\lra A^2(Y)
  $$
and the splitting
  $$
  A^2(X)=A^2(Y)\oplus J\; ,
  $$
where $J=I\cap A^2(X)$ inside $CH^2(X)$, if we fix a section
$A^2(Y)\to A^2(X)$ and identify $A^2(Y)$ with its image in $A^2(X)$.

For each index $j$ let $\tilde X_{p_j}$ be a resolution of
singularities of the surface $X_{p_j}$ ($\tilde X_{p_j}=X_{p_j}$ if
$X_{p_j}$ is non-singular). Note that for surfaces we have
resolution of singularities in any characteristic, so we need not to
make any restrictions on the ground field $k$ here. By Bertini's
theorem, for any irreducible component $Z$ of $\tilde X_{p_j}$ there
exists a smooth linear section $T$ of the Picard variety $P=\Pic
^0(Z)$ of $Z$. Applying Albanese functor to the embedding $T\hra P$
we get a surjective homomorphism $J=\Alb (C)\to \Alb (P)\cong P$
where $J$ is the Jacobian of the curve $C$. But any such
homomorphism can be induced by a divisor $D$ on the product $C\times
Z$ after tensoring with $\QQ $, see \cite{Scholl}. In other words,
$D_*:A^1(T)\to A^1(Z)$ is surjective. Since the push-forward of a
blowing up is surjective on Chow groups, \cite[6.7(b)]{Ful}, one can
easily construct a finite collection of non-singular curves
$T_1,\dots ,T_m$, and divisors $D_{ij}$ on $T_i\times X_{p_j}$, such
that the corresponding homomorphism $\sum _{i,j}(D_{ij})_*$ from
$\sum _iA^1(T_i)$ to $\sum _jA^1(X_{p_j})$ is surjective. Since
Neron-Severi group $NS_{\QQ }(Z)=CH^1(Z)/A^1(Z)$ is
finite-dimensional for each component $Z$ in $\tilde X_{p_j}$, the
complement of $\sum _jA^1(X_{p_j})$ in $\sum _jCH^1(X_{p_j})$ is a
finite-dimensional $\QQ $-vector space. It follows, that the kernel
$J$ of the localization homomorphism $A^2(X)\to A^2(Y)$ splits into
two $\QQ $-vector subspaces $V$ and $W$, where $W$ is covered by
$A^1$ of the above curves $T_i$ via the divisors $D_{ij}$, and the
second subspace $V$ is finite-dimensional. Then, by Lemma
\ref{uptofin}, and also taking into account the second definition of
representability of $A^2(X)$, we see that in order to prove that
$A^2(X)$ is representable we need only to show that there exists
another one finite collection of non-singular projective curves and
correspondences from them to $X$, such that the corresponding
homomorphism covers the complement $A^2(Y)$ of $J=V\oplus W$ inside
$A^2(X)$. In other words, in proving representability of $A^2(X)$ we
can cut out any finite number of fibers of the map $X\to C$ taking
into account the group $A^2(Y)$ only.

\bigskip

\noindent {\it Reduction 4}

\medskip

Extending $C$ more, if necessary, we may assume that $X_{\eta }$ has
a rational point over $\eta $. This rational point induces a section
of the map $f:Y\to U$ over some Zariski open subset $U'$ in $U$. If
we cut out a finite number of points in $U$, we omit a finite number
of non-singular fibers of the family $Y\to U$. Each non-singular
surface is regular, so that its group of divisors is just a
Neron-Severi group, which is finitely generated. Now again: since
representability of $A^2$ is modulo finite-dimensional subspaces in
$A^2$, see Lemma \ref{uptofin}, we may cut out any finite number of
points from $U$. It follows that, again, without loss of generality,
we may assume that the regular map $Y\to U$ has a section over $U$.

If $E$ is an image of that section, the self-intersection $E\cdot E$
is trivial on $Y$ because the codimension of $E$ in $X$ is equal to
two (or because of another one cutting of additional fibers of the
map $f$). Then it follows that the relative projectors
  $$
  \pi _0=E\times _UY\; ,\; \; \; \pi _4=Y\times _UE
  $$
are pairwise orthogonal, and hence
  $$
  \pi _2=\Delta _{Y/U}-\pi _0-\pi _4
  $$
is a relative second Murre projector for the whole family $f:Y\to
U$.

\bigskip

\noindent {\it Reduction 5}

\medskip

Let $M^2(Y/U)$ be the relative motive defined by the projector $\pi
_2$. Then we get the following standard decomposition:
  $$
  M(Y/U)=\uno _U\oplus M^2(Y/U)\oplus \Le _U^{\otimes 2}\; .
  $$
By the assumption of Theorem \ref{main}, the motive $M(X_{\bar \eta
})$ is finite-dimensional. If we look through the definition of
motivic finite-dimensionality from the geometrical viewpoint, it
means an existence of some algebraic cycles on $X_{\bar \eta
}^{\times 2N}\times \PR ^1$ providing rational triviality of the
wedge and symmetric powers of even and odd components of the
diagonal for $X_{\bar \eta }$. These algebraic cycles have their
common minimal field of definition, which is a finite extension of
$k(C)$. Extending $C$ by a finite extension again, we may assume,
without loss of generality, that the motive $M(Y_{\eta })$ is finite
dimensional over $\eta $ itself.

On the other hand, $\cm (\eta )$ is a colimit of the categories of
Chow motives over Zariski open subsets in $S$ by Lemma \ref{limit}.
It means that there exists $U$, such that the relative Chow motive
$M(Y/U)$ is finite-dimensional in the category $\cm (U)$.

\bigskip

\noindent {\it Main computations}

\medskip

Finite-dimensionality of $M(Y/U)$ implies, of course, that
$M^2(Y/U)$ is finite dimensional. But, actually, it is evenly finite
dimensional of dimension $b_2$. One can show that either using Lemma
\ref{limit} or by the following argument. Apriori one has
  $$
  M^2(Y/U)=A\oplus B\; ,
  $$
where $\wedge ^mA=0$ and $\Sym ^nB=0$. The base change functor
  $$
  \Xi :\cm (U)\lra \cm (\eta )\; .
  $$
is tensor, so it respects finite-dimensionality. In addition, $\Xi $
induces an isomorphism between $\End (\uno _U)=\QQ $ and $\End (\uno
_{\eta })=\QQ $. Since
  $$
  \Xi (M^2(Y/U))=M^2(Y_{\eta })
  $$
it follows that
  $$
  \Xi (B)=0
  $$
because $M^2(Y_{\eta })$ is evenly finite dimensional. Then
  $$
  B=0
  $$
by Lemma \ref{traces}. Thus, $M^2(Y/U)$ is evenly finite
dimensional. Using the the same arguments one can also show that it
can be annihilated by $\wedge ^{b_2+1}$.

\bigskip

The divisors $D_1,\dots ,D_{b_2}$ are defined over $\eta $ by
Reduction 2, and they generate the second cohomology group
$H^2(Y_{\bar \eta })$ via the cycle class map
  $$
  CH^1(X_{\eta })\lra H^2(Y_{\bar \eta })\; .
  $$
Since the motive $M(Y_{\eta })$ is finite dimensional, $H^1(Y_{\bar
\eta })=0$ and $H^2_{\tr }(Y_{\bar \eta })=0$, the second piece
$M^2(Y_{\eta })$ in the Murre decomposition of $M(Y_{\eta })$ can be
computed as follows:
  $$
  M^2(Y_{\eta })=\Le _{\eta }^{\oplus b_2}\; ,
  $$
see \cite[Theorem 2.14]{GP1}. Actually, these $b_2$ copies of $\Le
_{\eta }$ arise from the collection of divisor classes
  $$
  [D_1],\dots ,[D_{b_2}]
  $$
and their Poincar\'e dual classes
  $$
  [D'_1],\dots ,[D'_{b_2}]
  $$
on $Y_{\eta }$, loc.cit. In other words, if
  $$
  (\pi _2)_{\eta }=\Xi (\pi _2)
  $$
is a projector determining the middle motive $M^2(Y_{\eta })$, the
difference
  $$
  \varrho _{\eta }=(\pi _2)_{\eta }-\sum _{i=1}^{b_2}[D_i\times _{\eta }D_i']
  $$
is homologically trivial. Then
  $$
  \varrho _{\eta }^n=0
  $$
in the associative ring
  $$
  \End (M^2(Y_{\eta }))
  $$
for some $n$ by Kimura's nilpotency theorem (see Proposition 7.2 in
\cite{Kimura}).

Now let
  $$
  W_i\hspace{10mm} \hbox{and}\hspace{10mm} W'_i
  $$
be spreads of the above divisors $D_i$ and $D_i'$ over $U$
respectively (reduce $U$ some more, if necessary). By definition,
$W_i$ and $W_i'$ are algebraic cycles of codimension one on $Y$,
such that their pull-backs to the generic fiber are the divisors
$D_i$, and the same for $W'_i$. In other terms:
  $$
  \Xi [W_i]=D_i\; .
  $$
These $W_i$ and $W_i'$ are defined not uniquely, of course. The
cycles $W_i\times _UW_i'$ are in $\Corr _U^0(Y\times _UY)$, and we
set
  $$
  \varrho =\pi _2-
  \sum _{i=1}^{b_2}[W_i\times _UW_i']\; .
  $$
Then we have:
  $$
  \Xi (\varrho )=\varrho _{\eta }\; .
  $$

\bigskip

Let $\omega $ be any endomorphism of the motive $M^2(Y/U)$ and let
$\omega _{\eta }=\Xi (\omega )$. Then:
  $$
  \Xi (\tr (\omega \circ \varrho ))=\tr (\Xi (\omega \circ \varrho ))=
  \tr (\omega _{\eta }\circ \varrho _{\eta })=0
  $$
because $\varrho _{\eta }$ is homologically trivial. Here we use the
formula on page 116 in \cite{DeligneMilne}, i.e. the compatibility
of tensor functors with traces, once again. Since the functor $\Xi $
induces an isomorphism $\End (\uno _U)\cong \End (\uno _{\eta })$,
  $$
  \tr (\omega \circ \varrho )=0
  $$
for any $\omega $, i.e. $\varrho $ is numerically trivial.
Therefore,
  $$
  \varrho ^n=0
  $$
in $\End (M^2(Y/U))$ by Proposition \ref{AKnilp}.

Now let $\bar W_i$ be Zariski closure of $W_i$ in $X$ and let
  $$
  \theta _i=\Gamma _f^t\cdot [C\times \bar W_i]
  $$
be a cycle class of codimension $2$ in $C\times X$. Here the fibered
product $\times $ is assumed to be taken over $k$, and the cycle
class $\Gamma _f^t$ is a transpose of the graph of the map $f$, i.e.
the push-forward of the class $[X]$ with respect to the proper
morphism $\tau :X\to C\times X$ defined by the map $f:X\to C$ and
the identity $\id :X\to X$. Let
  $$
  (\theta _i)_*:CH^1(C)\lra CH^2(X)
  $$
be a homomorphism induced by the correspondence $\theta _i$. This
homomorphism can be computed also by another formula. Indeed, for
any cycle class $a$ from $CH^1(C)$ we have:
  $$
  f^*(a)=\tau ^*p_C^*(a)=\tau ^*p_C^*(a)\cdot [X]\; .
  $$
By the projection formula we get:
  $$
  \tau _*f^*(a)=
  p_C^*(a)\cdot \tau _*[X]=
  p_C^*(a)\cdot \Gamma _f^t\; .
  $$
Then we can compute:
  $$
  \begin{array}{rcl}
  (\theta _i)_*(a)
  &=&
  (p_X)_*(p_C^*(a)\cdot \theta _i) \\
  &=&
  (p_X)_*(p_C^*(a)\cdot \Gamma _f^t\cdot [C\times \bar W_i]) \\
  &=&
  (p_X)_*(\tau _*f^*(a)\cdot [C\times \bar W_i]) \\
  &=&
  (p_X)_*(\tau _*f^*(a)\cdot p_X^*[\bar W_i]) \\
  &=&
  (p_X)_*\tau _*f^*(a)\cdot [\bar W_i] \\
  &=&
  f^*(a)\cdot [\bar W_i]
  \end{array}
  $$
This gives us a possibility to define a homomorphism
  $$
  (\theta _i)_*:CH^1(U)\lra CH^2(Y)
  $$
in the non-compact case by the analogous formula:
  $$
  (\theta _i)_*(a)=f^*(a)\cdot [W_i]
  $$
for any cycle class $a$ in $CH^1(U)$, where $f$ is the map $f:Y\to
U$. Then, for each index $i$, we have a commutative square
  $$
  \diagram
  CH^1(C) \ar[dd]_-{} \ar[rr]^-{(\theta _i)_*} & & CH^2(X) \ar[dd]^-{} \\ \\
  CH^1(U) \ar[rr]^-{(\theta _i)_*} & & CH^2(Y)
  \enddiagram
  $$
In other words, the homomorphisms $(\theta _i)$ are compatible in
compact and non-compact cases.

Evidently, the homomorphisms $\theta _i$ respect algebraic
equivalence, so that we also have the corresponding commutative
square for groups $A^*(-)$:
  $$
  \diagram
  A^1(C) \ar[dd]_-{} \ar[rr]^-{(\theta _i)_*} & & A^2(X) \ar[dd]^-{} \\ \\
  A^1(U) \ar[rr]^-{(\theta _i)_*} & & A^2(Y)
  \enddiagram
  $$
A resulting commutative diagram we need is then as follows:
  $$
  \diagram
  \oplus _{i=1}^{b_2}A^1(C) \ar[dd]_-{} \ar[rr]^-{\theta _*} & &
  A^2(X) \ar[dd]^-{} \\ \\
  \oplus _{i=1}^{b_2}A^1(U) \ar[rr]^-{\theta _*} & & A^2(Y)
  \enddiagram
  $$
Both homomorphisms $\theta _*$ here are defined by summing of
homomorphisms $(\theta _i)_*$.

An elementary diagram chasing shows now that, if the bottom
homomorphism $\theta _*$ is surjective, the image of the top
homomorphism $\theta _*$ coincides with $A^2(X)$ modulo a
finite-dimensional subspace $J$ in $A^2(X)$ appeared in Reduction 3.
Let us recall that, as soon as we fix a section $A^2(Y)\to A^2(X)$
we have that $A^2(X)=A^2(Y)\oplus J$, $J=W\oplus V$, $W$ is covered
by $A^1$ of curves and $V$ is finite-dimensional. Then, according to
what was shown in Reduction 3, in order to prove representability of
$A^2(X)$ we need only to show that the bottom $\theta _*$ is onto.

Let $y\in CH^i(Y)$. For any correspondence $c\in \Corr _U^j(Y,Y)$
let, as usual,
  $$
  c_*(y)={p_2}_*(p_1^*(y)\cdot c)\in CH^{i+j}(Y)
  $$
be the action of the correspondence $c$ on $y$ in the relative
sense, i.e. $p_1$ and $p_2$ are two projections from $Y\times _UY$
onto $Y$. In particular, one has a decomposition
  $$
  \begin{array}{rcl}
  y
  &=&
  {\Delta _{Y/U}}_*(y) \\
  &=&
  (\pi _0)_*(y)+(\pi _2)_*(y)+(\pi _4)_*(y)
  \end{array}
  $$
in $CH^2(Y)$. Now we have:
  $$
  \begin{array}{rcl}
  (\pi _0)_*(y)
  &=&
  {p_2}_*(p_1^*(y)\cdot \pi _0) \\
  &=&
  {p_2}_*(p_1^*(y)\cdot p_1^*[E]) \\
  &=&
  {p_2}_*p_1^*(y\cdot [E]) \\
  &=&
  f^*f_*(y\cdot [E])
  \end{array}
  $$
From now on we assume that $y$ is of codimension two. Then, as $y$
and $[E]$ are both of codimension two in a three-dimensional
variety, we have that $y\cdot [E]=0$, whence
  $$
  (\pi _0)_*(y)=0\; .
  $$
Assume, furthermore, that $y$ is algebraically trivial, i.e. $y\in
A^2(Y)$. In that case $f_*(y)=0$. Then we compute:
  $$
  \begin{array}{rcl}
  (\pi _4)_*(y)
  &=&
  {p_2}_*(p_1^*(y)\cdot \pi _4) \\
  &=&
  {p_2}_*((y\times _U[Y])\cdot ([Y]\times _U[E])) \\
  &=&
  {p_2}_*(y\times _U[E]) \\
  &=&
  f_*(y)\times _U[E] \\
  &=&
  0
  \end{array}
  $$
(here and below we use nice properties of intersections of relative
cycles over a non-singular one-dimensional base, see Chapter 2.2 in
\cite{Ful}). As a result we have:
  $$
  y=(\pi _2)_*(y)\; .
  $$
On the other hand,
  $$
  \pi _2=\sum _{i=1}^{b_2}[W_i\times _UW_i']+\varrho \; ,
  $$
whence we get:
  $$
  y=(\pi _2)_*(y)=
  \sum _{i=1}^{b_2}[W_i\times _UW_i']_*(y)
  +\varrho _*(y)\; .
  $$
Let us emphasize once more that the correspondences here act as
relative correspondences.

Now let
  $$
  v_1=-\sum _{i=1}^{b_2}[W_i\times _UW_i']_*(y)\; ,
  $$
so that
  $$
  \varrho _*(y)=y+v_1\; .
  $$
Write $y$ as a class of a linear combination
  $$
  \sum _jn_jZ_j\; ,
  $$
where $Z_j$ are integral curves on $Y$. For any $i$ and $j$ one has:
  $$
  \begin{array}{rcl}
  [W_i\times _UW_i']_*[Z_j]
  &=&
  {p_2}_*([Z_j\times _UY]\cdot [W_i\times _UW_i']) \\
  &=&
  {p_2}_*(([Z_j]\cdot [W_i])\times _U([Y]\cdot [W_i'])) \\
  &=&
  {p_2}_*(([Z_j]\cdot [W_i])\times _U[W_i'])
  \end{array}
  $$
By linearity:
  $$
  [W_i\times _UW_i']_*(y)={p_2}_*((y\cdot [W_i])\times _U[W_i'])
  $$
Since $y$ is of codimension two and $W_i$ is of codimension one in
the non-singular threefold $Y$, the intersection $y\cdot [W_i]$ is
zero-dimensional cycle class on $Y$. Let
  $$
  a_i=f_*(y\cdot [W_i])
  $$
be its push-forward to $U$ with respect to the proper map $f:Y\to
U$. Using proper-flat base change and the projection formula we
compute:
  $$
  \begin{array}{rcl}
  [W_i\times _UW_i']_*(y)
  &=&
  {p_2}_*((y\cdot [W_i])\times _U[W_i']) \\
  &=&
  {p_2}_*(p_1^*(y\cdot [W_i])\cdot p_2^*([W_i'])) \\
  &=&
  {p_2}_*p_1^*(y\cdot [W_i])\cdot [W_i'] \\
  &=&
  f^*f_*(y\cdot [W_i])\cdot [W_i'] \\
  &=&
  f^*(a_i)\cdot [W_i'] \\
  &=&
  (\theta _i)_*(a_i)\; .
  \end{array}
  $$
Moreover, as $y$ is in $A^2(Y)$, it follows that each $a_i$ is in
$A^1(U)$. Then we get:
  $$
  -v_1=\sum _{i=1}^{b_2}[W_i\times _UW_i']_*(y)=
  \sum _{i=1}^{b_2}(\theta _i)_*(a_i)=\theta _*(e_1)\; ,
  $$
where
  $$
  c_1=(a_1,\dots ,a_{b_2})
  $$
is in $\oplus _{i=1}^{b_2}A^1(U)$.

Applying $\varrho _*$ to the both sides of the equality
  $$
  \varrho _*(y)=\theta _*(c_1)+y
  $$
we get:
  $$
  \varrho _*^2(y)=\theta _*(c_2)+2y
  $$
for some $c_2$ in $\oplus _{i=1}^{b_2}A^1(U)$, and so forth. After
$n$ steps we will get:
  $$
  \varrho _*^n(y)=\theta _*(c_n)+ny\; ,
  $$
where $e_n$ is in $\oplus _{i=1}^{b_2}A^1(U)$. But we know that
$\varrho $ is a nilpotent correspondence in $\End (M(Y/U))$, so that
  $$
  \varrho ^n_*(y)=0
  $$
for big enough $n$. It follows that
  $$
  y=-\frac{1}{n}\cdot \theta _*(c_n)
  $$
is in the image of the homomorphism $\theta _*$, and we are done.

\bigskip

{\it Acknowledgements.} The author wishes to thank Joseph Ayoub,
Spencer Bloch, H\'el\`ene Esnault, Aleksandr Pukhlikov and Andrei
Suslin for useful discussions on the theme of that paper. The work
was supported by NSF Grant DMS-0111298. I am grateful to the
Institute for Advanced Study at Princeton for the support and
hospitality in 2005-2006, and to JSPS, Tokyo University and
Hiroshima University for funding of my stay in Japan in June - July
2005.

\bigskip

\begin{small}

\end{small}

\vspace{8mm}

\begin{small}

{\sc Department of Mathematical Sciences, University of Liverpool,
Peach Street, Liverpool L69 7ZL, England, UK}

\medskip

\begin{footnotesize}

{\it E-mail address}: {\tt vladimir.guletskii@liverpool.ac.uk}

\end{footnotesize}

\end{small}

\end{document}